# WEAK CONVERGENCE OF POSITIVE SELF-SIMILAR MARKOV PROCESSES AND OVERSHOOTS OF LÉVY PROCESSES

By M. E. Caballero and L. Chaumont

*Universidad Nacional Autónoma de México and Université Pierre et Marie Curie*

Using Lamperti's relationship between Lévy processes and positive self-similar Markov processes (pssMp), we study the weak convergence of the law $\mathbb{P}_x$ of a pssMp starting at $x > 0$, in the Skorohod space of càdlàg paths, when $x$ tends to 0. To do so, we first give conditions which allow us to construct a càdlàg Markov process $X^{(0)}$, starting from 0, which stays positive and verifies the scaling property. Then we establish necessary and sufficient conditions for the laws $\mathbb{P}_x$ to converge weakly to the law of $X^{(0)}$ as $x$ goes to 0. In particular, this answers a question raised by Lamperti [*Z. Wahrsch. Verw. Gebiete* **22** (1972) 205–225] about the Feller property for pssMp at $x = 0$.

**1. Introduction.** An $\mathbb{R}_+$-valued self-similar Markov process $X$, under the family of probabilities $(\mathbb{P}_x, x \geq 0)$ is a càdlàg Markov process which fulfills a scaling property, that is, there exists a constant $\alpha > 0$ such that

(1.1)     the law of $(kX_{k^{-\alpha}t}, t \geq 0)$ under $\mathbb{P}_x$ is $\mathbb{P}_{kx}$     for all $k > 0$.

Self-similar Markov processes are involved in various parts of probability theory, such as branching processes and fragmentation theory. They also arise as limits of re-scaled Markov processes. Their properties have been deeply studied by the early 1960s, especially through Lamperti's work on one-dimensional branching processes. In this paper we focus on positive self-similar Markov processes to which we will refer as pssMp. Some particularly well-known examples which are discussed in Section 4 are Bessel processes, stable subordinators, or more generally, stable Lévy processes conditioned to stay positive.

---

From a straightforward argument using (1.1), it is easily seen that the family of measures $(\mathbb{P}_x)$ defined on Skorohod's space of càdlàg functions is weakly continuous in $x$ on the open half-line $(0, \infty)$. The question of the existence of a weak limit when $x$ goes toward $0$ is much less obvious. This is the main object of our paper.

Let $\rho \overset{(\mathrm{def})}{=} \inf\{t : X_t = 0\}$, then it appears from Lamperti's study [11] that any positive self-similar Markov process $X$ is such that either $\rho < \infty$ and $X_{\rho-} = 0$ a.s., or $\rho < \infty$ and $X_{\rho-} > 0$ a.s., or $\rho = +\infty$ a.s., and this trichotomy does not depend on the starting point $x > 0$. In the two first cases, the Markov property implies that the process $(X_{\rho+t}, t \geq 0)$ is independent of the process $(X_t, t \leq \rho)$ and its law does not depend on $x$. Moreover, the scaling property implies that $\lim_{x\to 0} \rho = 0$ a.s., hence, this shows that the family of measures $(\mathbb{P}_x)$ converges weakly, as $x$ goes to $0$, toward the law of the process $(X_{\rho+t}, t \geq 0)$. Let us mention that in these cases Rivero [12, 13] studied the different ways to construct an entrance law for the process $(X_{\rho+t}, t \geq 0)$. So we shall focus on the last case, that is, when $\rho = +\infty$, a.s. Then the state space of $(X, \mathbb{P}_x)$, $x > 0$, is $(0, \infty)$ and $0$ is a boundary point. It is then natural to wonder if the semigroup of $(X, \mathbb{P}_x)$ can be extended to the nonnegative half-line $[0, \infty)$; in other words, can an entrance law at $0$ for $(X, \mathbb{P}_x)$ be defined? This problem has first been raised by Lamperti [11], who observed in the proof of his Theorem 2.1 that "... *the Feller property may fail at* $x = 0$." It has been partially solved by Bertoin and Caballero [2] and Bertoin and Yor [4], who gave sufficient conditions for the weak convergence of $\mathbb{P}_x$ to hold when $x$ tends to $0$, in the sense of finite-dimensional distributions. In this paper, we characterize the self-similar families of laws $(\mathbb{P}_x)$ which converges weakly as $x$ tends to $0$, on the Skorohod space of càdlàg functions. We also describe their limit law $\mathbb{P}_0$ by constructing the paths of a self-similar Markov process whose law is $\mathbb{P}_0$.

A crucial point in our arguments is the famous Lamperti representation of self-similar $\mathbb{R}_+$-valued processes. This transformation enables us to construct the paths of any such self-similar process $X$ from those of a Lévy process. More precisely, Lamperti [11] found the representation

$$(1.2) \qquad X_t = x\exp\xi_{\tau(tx^{-\alpha})}, \qquad t \geq 0,$$

under $\mathbb{P}_x$, for $x > 0$, where

$$\tau_t = \inf\{s : A_s > t\}, \qquad A_s = \int_0^s \exp\alpha\xi_u\, du,$$

and where $\xi$ is either a real Lévy process such that $\lim_t \xi_t = -\infty$, if $\rho < +\infty$ and $X_{\rho-} = 0$, or $\xi$ is a Lévy process killed at an independent exponential time if $\rho < +\infty$ and $X_{\rho-} > 0$, or $\xi$ is a Lévy process such that $\limsup_t \xi_t = +\infty$, if $\rho = +\infty$. Note that, for $t < A_\infty$, we have the equality $\tau_t = \int_0^t \frac{ds}{X_s^\alpha}$, so that (1.2) is invertible and yields a one to one relation between the class of pssMp and the one of Lévy processes.

The conditions of weak convergence of the family $(\mathbb{P}_x)$ and the construction of the limit process which are presented hereafter will naturally bear on the features of the underlying Lévy process in the Lamperti transformation. The rest of the paper is organized in three sections. In Section 2 the main results are stated and discussed. The proofs are given in Section 3 and some examples are presented in Section 4.

## 2. Main results.

2.1. *Preliminaries on overshoots.* On the Skorohod space $\mathcal{D}$ of càdlàg trajectories we consider a reference probability measure $P$ under which $\xi$ will always denote a real Lévy process such that $\xi_0 = 0$. We now give some preliminaries on overshoots of Lévy processes, that is, $(\xi_{T_z} - z, z \geq 0)$ with $T_z = \inf\{t : \xi_t \geq z\}$. The condition of weak convergence of the overshoot process $(\xi_{T_z} - z, z \geq 0)$ when $z$ tends to $+\infty$ will appear naturally in our main results (Theorems 1 and 2) to be necessary for the probability measures $(\mathbb{P}_x)$ to converge weakly as $x \to 0$. The asymptotic behavior of the overshoot process of Lévy processes has strongly been studied by Doney and Maller [10]; let us briefly recall one of their main results. Let $\Pi$ be the Lévy measure of $\xi$, that is, the measure satisfying $\int_{(-\infty,\infty)}(x^2 \wedge 1) \times \Pi(dx) < \infty$ and such that the characteristic exponent $\Psi$ [defined by $E(e^{iu\xi_t}) = e^{-t\Psi(u)}$, $t \geq 0$] is given, for some $b \geq 0$ and $a \in \mathbb{R}$, by

$$\Psi(u) = iau + \tfrac{1}{2}b^2 u^2 + \int_{(-\infty,\infty)}(1 - e^{iux} + iux\mathbb{1}_{\{|x| \leq 1\}})\Pi(dx), \qquad u \in \mathbb{R}.$$

Define, for $x \geq 0$,

$$\overline{\Pi}^+(x) = \Pi((x,\infty)), \qquad \overline{\Pi}^-(x) = \Pi((-\infty,-x))$$

and

$$J = \int_{[1,\infty)} \frac{x\overline{\Pi}^+(x)\,dx}{1 + \int_0^x dy \int_y^\infty \overline{\Pi}^-(z)\,dz}.$$

Then, according to Doney and Maller [10], a necessary and sufficient condition for the overshoot $\xi_{T_z} - z$ to converge weakly toward the law of a finite r.v. as $z$ goes to $+\infty$ is

(H)    $\xi$ is not arithmetic and $\begin{cases} \text{either } 0 < E(\xi_1) \leq E(|\xi_1|) < \infty, \\ \text{or } E(|\xi_1|) < \infty, \ E(\xi_1) = 0 \text{ and } J < \infty. \end{cases}$

(Actually, the case where $\xi$ has no negative jumps and the case where $\xi$ is arithmetic are not considered in [10]. This makes our expression of the integral $J$ slightly different from the one which is given in [10].) The discrete time case has been treated in [8].

As noted in [10], the above condition may also be expressed in terms of the upward ladder height process, say, $\sigma$, associated with $\xi$. (We refer to [1], Chapter VII for a definition of the ladder process.) Indeed, we easily see that, for each level $x$, the overshoots of the processes $\xi$ and $\sigma$ across $x$ are the same, hence, (H) is equivalent to say that $\sigma$ is not arithmetic and $E(\sigma_1) < \infty$. Let us also mention that the latter result is proved for any subordinator in [3] where the authors also give the explicit law of the limit of the overshoots. In the rest of this paper, under hypothesis (H), $\theta$ will denote an r.v. whose law is the weak limit of the overshoots of $\xi$, that is,

$$(2.1) \qquad \xi_{T_z} - z \xrightarrow[z \to \infty]{(w)} \theta.$$

Then from [3] the law of $\theta$ is given by

$$(2.2) \qquad \theta \stackrel{(d)}{=} \mathcal{U} Z,$$

where $\mathcal{U}$ and $Z$ are independent r.v.'s, $U$ is uniformly distributed over $[0, 1]$ and, $\nu$ being the Lévy measure of the subordinator $\sigma$, the law of $Z$ is given by

$$\mathbb{P}(Z > t) = \mathbb{E}(\sigma_1)^{-1} \int_{(t, \infty)} s \nu(ds), \qquad t \geq 0.$$

In the sequel we will need a multidimensional version of the convergence which is stated in (2.1). As the next proposition shows, it is equivalent to the stationarity of the process $(\xi_{T_{z-\theta}} - (z - \theta), z \geq 0)$.

PROPOSITION 1. *If the overshoot process $(\xi_{T_z} - z, z \geq 0)$ converges weakly at infinity toward the law of an a.s. finite r.v. $\theta$, then the multidimensional overshoot converges weakly. More precisely, for every increasing sequence of nonnegative reals $0 \leq z_1 \leq \cdots \leq z_k$,*

$$(\xi_{T_{z+z_1}} - (z + z_1), \ldots, \xi_{T_{z+z_k}} - (z + z_k))$$

$$\xrightarrow[z \to \infty]{(w)} (\xi_{T_{z_1-\theta}} - (z_1 - \theta), \ldots, \xi_{T_{z_k-\theta}} - (z_k - \theta)),$$

*where $\theta$ is independent of $\xi$. In particular, $(\xi_{T_{z-\theta}} - (z - \theta), z \geq 0)$ is a strictly stationary Markov process, that is, for any nonnegative real $r$, $(\xi_{T_{z+r-\theta}} - (z + r - \theta), z \geq 0) \stackrel{(d)}{=} (\xi_{T_{z-\theta}} - (z - \theta), z \geq 0)$.*

*Conversely, if there exists an a.s. finite r.v. $\theta$, independent of $\xi$ such that $(\xi_{T_{z-\theta}} - (z - \theta), z \geq 0)$ is strictly stationary, then the overshoot process $(\xi_{T_z} - z, z \geq 0)$ converges weakly at infinity toward the law of $\theta$.*

Then an important point in the construction of a weak limit for $(\mathbb{P}_x)$ is the following definition of a random sequence of $(\mathcal{D} \times \mathbb{R}_+)^{\mathbb{N}}$.

PROPOSITION 2. *Let $x_1 \geq x_2 \geq \cdots > 0$ be an infinite decreasing sequence of positive real numbers which converges toward $0$. Under condition* (H), *there exists a random sequence $(\theta_n, \xi^{(n)})$ of $(\mathbb{R}_+ \times \mathcal{D})^{\mathbb{N}}$ such that, for each $n$, $\theta_n$ and $\xi^{(n)}$ are independent, $\theta_n \stackrel{(d)}{=} \theta$, $\xi^{(n)} \stackrel{(d)}{=} \xi$, and for any $i, j$ such that $1 \leq i \leq j$,*

$$(2.3) \quad \xi^{(i)} \stackrel{(a.s.)}{=} (\xi^{(j)}(T^{(j)}_{\log(x_i e^{-\theta_j}/x_j)} + t) - \xi^{(j)}(T^{(j)}_{\log(x_i e^{-\theta_j}/x_j)}), t \geq 0),$$

$$(2.4) \quad \theta_i \stackrel{(a.s.)}{=} \xi^{(j)}(T^{(j)}_{\log(x_i e^{-\theta_j}/x_j)}) - \log(x_i e^{-\theta_j}/x_j),$$

*where $T^{(j)}_z = \inf\{t : \xi^{(j)}_t \geq z\}$, $z \in \mathbb{R}$. The above statement determines the law of the sequence $(\xi^{(n)}, \theta_n)$. Furthermore, for any $n$, $\xi^{(n)}$ is independent of $(\theta_k)_{k \geq n}$ and $(\theta_n)$ is a Markov chain.*

A particularity of the above definition is the "backward" inductive construction of the law of $(\theta_n, \xi^{(n)})$. More specifically, we see from (2.3) and (2.4) that, for any $n$, the couples $(\theta_{n-1}, \xi^{(n-1)}), (\theta_{n-2}, \xi^{(n-2)}), \ldots, (\theta_1, \xi^{(1)})$ are functionals of $(\theta_n, \xi^{(n)})$.

As will be seen in the next subsection, for each $n$, $(\theta_n, \xi^{(n)})$ is used to construct our limiting process between its first passage time above $x_n$ and its first passage time above $x_{n-1}$. Then the result of Proposition 2 will allow us to describe the joint law of the values of this process at its first passage times above $x_n$, $n \geq 1$.

We emphasize that if $\xi$ has no positive jumps, then $\theta_i = 0$ a.s. and the results of this paper, as well as their proofs, are much simpler. Note also that, in any other case, it may happen that the event

$$\{T^{(j)}_{\log(x_i e^{-\theta_j}/x_j)} = 0\} = \{\log(x_i e^{-\theta_j}/x_j) \leq 0\}$$

has positive probability. On this event, we have $x_i e^{\theta_i} = x_j e^{\theta_j}$ and $\xi^{(i)} = \xi^{(j)}$.

2.2. *The process issued from* 0. We recall that $P$ is our reference probability measure and we denote by $X$ the canonical process on $\mathcal{D}$.

Consider any self-similar Markov process $(X, \mathbb{P}_x)$, $x > 0$, as defined in (1.2), with $\rho = \inf\{t : X_t = 0\} = +\infty$, $\mathbb{P}_x$-a.s. As we recalled in the Introduction, there exists a unique Lévy process $\xi$ with $\limsup_{t \to +\infty} \xi_t = +\infty$, a.s. and such that the representation (1.2) holds. Observe that since the time change $t \mapsto \tau(tx^{-\alpha})$ is continuous with $\lim_{t \to +\infty} \tau(tx^{-\alpha}) = +\infty$, a.s., and since $x \mapsto e^x$ is increasing, we have $\mathbb{P}_x$-a.s., $\limsup_{t \to +\infty} X_t = +\infty$.

Note that from (1.1), the process $(X^\alpha, \mathbb{P}_x)$, $x > 0$, is a pssMp whose scaling coefficient is equal to 1. Moreover, the power function is a continuous functional of the càdlàd trajectories, hence, we do not lose any generality in our convergence type results by assuming that $\alpha = 1$.

In the next theorems we first give the construction of a Markov process $X^{(0)}$ which starts from 0 continuously, stays positive and fulfills the scaling property with index $\alpha = 1$. We shall then see that the law of the process $X^{(0)}$ is the weak limit on $\mathcal{D}$ of $X$, as $x$ tends to 0.

From the sequence $(\theta_n, \xi^{(n)})$ defined in Proposition 2, we first introduce the sequence of processes

$$(2.5) \qquad X_t^{(\bar{x}_n)} = \bar{x}_n \exp \xi_{\tau^{(n)}(t/\bar{x}_n)}^{(n)}, \qquad t \geq 0, n \geq 1,$$

where $\bar{x}_n = x_n e^{\theta_n}$ and with the natural definition $\tau_t^{(n)} \stackrel{(\text{def})}{=} \inf\{s : \int_0^s \exp \xi_u^{(n)} du > t\}$. Let also

$$(2.6) \qquad S^{(n-1)} = \inf\{t \geq 0 : X_t^{(\bar{x}_n)} \geq x_{n-1}\}, \qquad n \geq 2,$$

which is a.s. finite from our assumptions. Under an additional condition which ensures that $\Sigma_n = \sum_{k=n}^{+\infty} S^{(k)} < +\infty$, a.s., in the next theorem we define the process $X^{(0)}$ on the positive halfline as the concatenation of the processes $(X_{t-\Sigma_n}^{(\bar{x}_n)}, \Sigma_n \leq t \leq \Sigma_{n-1})$, $n \geq 2$, and $(X_{t-\Sigma_1}^{(\bar{x}_1)}, t \geq \Sigma_1)$.

THEOREM 1. *Assume that* (H) *holds and that* $E(\log^+ \int_0^{T_1} \exp \xi_s \, ds) < \infty$. *Let* $\Sigma_n = \sum_{k=n}^\infty S^{(k)}$, *then, for any* $n$, $0 < \Sigma_n < \infty$, *a.s., so that the following construction*

$$(2.7) \qquad X_t^{(0)} = \begin{cases} X_{t-\Sigma_1}^{(\bar{x}_1)}, & t \in [\Sigma_1, \infty), \\ X_{t-\Sigma_2}^{(\bar{x}_2)}, & t \in [\Sigma_2, \Sigma_1), \\ \vdots \\ X_{t-\Sigma_n}^{(\bar{x}_n)}, & t \in [\Sigma_n, \Sigma_{n-1}), \\ \vdots \end{cases} \qquad X_0^{(0)} = 0$$

*makes sense and it defines a càdlàg stochastic process on the real halfline* $[0, \infty)$ *with the following properties:*

(i) *The law of* $X^{(0)}$ *does not depend on the sequence* $(x_n)$.

(ii) *The paths of the process* $X^{(0)}$ *are such that* $\limsup_t X_t^{(0)} = +\infty$, *a.s. and* $X_t^{(0)} > 0$, *a.s. for any* $t > 0$.

(iii) *The process* $X^{(0)}$ *satisfies the scaling property, that is, for any* $k > 0$,

$$(2.8) \qquad (k X_{k^{-1}t}^{(0)}, t \geq 0) \text{ has the same law as } X^{(0)}.$$

(iv) *The process $X^{(0)}$ is strong Markov and has the same semigroup as $(X, \mathbb{P}_x)$ for $x > 0$.*

*We will denote by $\mathbb{P}_0$ the law of this process.*

Note that $\Sigma_n$ is the first passage time above $x_n$ by the process $X^{(0)}$, that is,

$$\Sigma_n = \inf\{t : X_t^{(0)} \geq x_n\},$$

and that the r.v. $\bar{x}_n = x_n e^{\theta_n}$ represents the value of this process at time $\Sigma_n$. In particular, when $\xi$ has no positive jumps, the above construction of $X^{(0)}$ says that $(X_{\Sigma_n+t}^{(0)}, t \geq 0)$ has law $\mathbb{P}_{x_n}$.

For the description of the law $\mathbb{P}_0$ of the process $X^{(0)}$ to be complete, it is worth describing its entrance law. As we will see in Theorem 2, $\mathbb{P}_0$ is the weak limit of $(\mathbb{P}_x)$ as $x$ tends to 0. So, in the case where $0 < E(\xi_1) \leq E(|\xi_1|) < \infty$, (i.e., $\xi$ has positive drift) the entrance law under $\mathbb{P}_0$ has been computed by Bertoin, Caballero [2] and Bertoin, Yor [4] and can be expressed as follows: for every $t > 0$ and for every measurable function $f : \mathbb{R}_+ \to \mathbb{R}_+$,

$$(2.9) \qquad \mathbb{E}_0(f(X_t)) = \frac{1}{m} E(I^{-1} f(tI^{-1})),$$

where $m = E(\xi_1)$ and $I = \int_0^\infty \exp(-\xi_s)\, ds$. When $E(\xi_1) = 0$, we have no explicit computation for the entrance law of $X^{(0)}$ in terms of the underlying Lévy process. However, in the following result, we see that it can be obtained as the weak limit of the entrance law for the positive drift case, when the drift tends toward 0.

PROPOSITION 3. *Under conditions of Theorem* 1, *suppose that $E(\xi_1) > 0$, then the entrance law under $\mathbb{P}_0$ is given by* (2.9). *Suppose that $E(\xi_1) = 0$, then for any bounded and continuous function $f$,*

$$(2.10) \qquad \mathbb{E}_0(f(X_t)) = \lim_{\lambda \to 0} \frac{1}{\lambda} E(I_\lambda^{-1} f(tI_\lambda^{-1})),$$

*where $I_\lambda = \int_0^\infty \exp(-\xi_s - \lambda s)\, ds$.*

The next theorem gives necessary and sufficient conditions bearing upon $\xi$, for a family of pssMp $(X, \mathbb{P}_x)$ described by (1.2) to converge weakly in the Skorohod space $\mathcal{D}$ of càdlàg trajectories, as $x$ goes toward 0. Recall that the first passage time process of $\xi$ is defined by $T_z = \inf\{t : \xi_t \geq z\}$, $z \in \mathbb{R}$. In the sequel, by degenerate probability measure on $\mathcal{D}$, we mean the law of a constant process which is finite or infinite.

THEOREM 2. *Let $(X, \mathbb{P}_x)$, $x > 0$, be defined as in (1.2) and such that $\rho = +\infty$ (or, equivalently, $\limsup_t \xi_t = +\infty$).*

*The family of probability measures $(\mathbb{P}_x)$ converges weakly in $\mathcal{D}$, as $x \to 0$, toward a nondegenerate probability measure if and only if the overshoot process $(\xi_{T_z} - z, z \geq 0)$ converges weakly toward the law of a finite r.v. as $z \to +\infty$ [i.e., (H) holds] and $E(\log^+ \int_0^{T_1} \exp \xi_s \, ds) < +\infty$. Under these conditions, the limiting law of $(\mathbb{P}_x)$ is $\mathbb{P}_0$.*

A consequence of Theorem 2 is that, under these conditions, the semi-group $P_t f(x) \overset{(\mathrm{def})}{=} \mathbb{E}_x(f(X_t))$ is Fellerian on the space $\mathcal{C}_0(\mathbb{R}_+)$ of continuous functions $f : \mathbb{R}_+ \to \mathbb{R}$ with $\lim_{x \to +\infty} f(x) = 0$. It will appear along the proofs in the next section that when (H) holds but $E(\log^+ \int_0^{T_1} \exp \xi_s \, ds) = \infty$, the process $(X, \mathbb{P}_x)$ actually converges weakly toward the process which is identically nought (see the remark at the end of Section 3). It is shown in [4] that when $E(\xi_1)$ exists but is infinite, the convergence of $(X, \mathbb{P}_x)$ toward a nondegenerate process in the sense of finite-dimensional distributions does not hold. The problem of finding necessary and sufficient conditions for the convergence in the sense of finite-dimensional distributions remains open. More precisely, we do not know what happens when $E(\xi_1)$ does not exist.

Before we prove the above results, let us discuss the condition

$$(2.11) \qquad E\left(\log^+ \int_0^{T_1} \exp \xi_s \, ds\right) < \infty$$

which is involved in Theorems 1 and 2. This condition is rather weak in the sense that it is satisfied by a very large class of Lévy processes. It contains at least the cases which are described hereafter.

First of all, if $0 < E(\xi_1) \leq E(|\xi_1|) < \infty$, then it is well known that $E(T_1) < \infty$ and, hence, condition (2.11) is satisfied.

Condition (2.11) also holds when $\xi$ satisfies the Spitzer condition, that is,

$$\frac{1}{t} \int_0^t \mathbb{1}_{\{\xi_s \geq 0\}} \, ds \to \rho \in (0, 1).$$

Indeed, in that case it is known (see [1], Theorem VI.18) that $E(T_1^\gamma) < +\infty$, for $0 < \gamma < \rho$. Hence, $\int_0^{T_1} \exp \xi_s \, ds (\leq e T_1)$ has a finite logarithmic moment.

Another instance when this condition is satisfied is when $\xi$ has an exponential moment, that is, $E(\exp \xi_1) < \infty$. Under this condition, there exists $\alpha > 0$ such that the process $\exp(\xi_{t \wedge T_1} - \alpha(t \wedge T_1))$ is a martingale and we have $E(\exp(\xi_{T_1} - \alpha T_1)) = 1$. (These conditions are satisfied, e.g., whenever $\xi$ has no positive jumps.) Then we can write

$$E(\mathbb{1}_{\{t < T_1\}} \exp \xi_t) = E(\mathbb{1}_{\{t < T_1\}} \exp[\xi_{T_1} - \alpha(T_1 - t)])$$

so that, integrating with respect to $t$, we obtain the stronger condition:

$$(2.12) \qquad E\left(\int_0^{T_1} \exp \xi_t \, dt\right) = \frac{1}{\alpha}(E(\exp \xi_{T_1}) - 1) < \infty.$$

It will appear in the proof of Lemma 2(iii) that condition (2.11) is equivalent to the following:

$$(2.13) \qquad E\left(\log^+ \int_0^{T_x} \exp \xi_s \, ds\right) < \infty \qquad \text{for any } x > 0.$$

To end this section, we emphasize that we do not have any example of a Lévy process $\xi$ for which (2.11) is not satisfied and this problem remains open.

## 3. Proofs.

PROOF OF PROPOSITION 1. Put $\theta^{(z)} = \xi_{T_z} - z$ and $\bar\xi = (\xi_{T_z+t} - \xi_{T_z}, t \geq 0)$ and note the path by path identity:

$$(\xi_{T_{z+z_1}} - (z + z_1), \ldots, \xi_{T_{z+z_n}} - (z + z_n))$$
$$= (\bar\xi(T_{z_1 - \theta^{(z)}}) - (z_1 - \theta^{(z)}), \ldots, \bar\xi(T_{z_n - \theta^{(z)}}) - (z_n - \theta^{(z)})).$$

As mentioned above, under hypothesis (H), $\theta^{(z)}$ converges weakly toward $\theta$. Moreover, the process $\bar\xi$ is distributed as $\xi$, so the weak convergence follows from the independence between $\bar\xi$ and $\theta^{(z)}$. The process $(\xi_{T_{z-\theta}} - (z - \theta), z \geq 0)$ is obviously càdlàg and its Markov property follows from general properties of Lévy processes and from the independence between $\theta$ and $\xi$.

Now suppose that there exists an a.s. finite r.v. $\theta$ which is independent of $\xi$ and such that the process $(\xi_{T_{z-\theta}} - (z - \theta), z \geq 0)$ is strictly stationary. Let $\sigma$ be the upward ladder height process associated to $\xi$, then as we already observed in Section 2.1, the overshoot of $\sigma$ and the overshoot of $\xi$ at any level $z$ are equal. Define the first passage process $\nu_z = \inf\{t : \sigma_t \geq z\}$, then $(\sigma_{\nu_{z-\theta}} - (z - \theta), z \geq 0)$ is strictly stationary. But it is shown in [10], at the end of the proof of Lemma 7, that if $E(\sigma_1) = +\infty$, then $\sigma_{\nu_z} - z$ tends to $+\infty$ in probability as $z$ goes to $+\infty$. This is in contradiction with the stationarity of $(\sigma_{\nu_{z-\theta}} - (z - \theta), z \geq 0)$, hence, $E(\sigma_1) < +\infty$, which is equivalent to (H), as mentioned in Section 2.1. □

We emphasize that the stationary property of the process $(\xi_{T_{z-\theta}} - (z - \theta), z \geq 0)$ is a crucial point in the proofs of our results, so we will often make use of Proposition 1 in the sequel.

PROOF OF PROPOSITION 2. Let $x_1 \geq x_2 \geq \cdots > 0$ be an infinite decreasing sequence of positive real numbers which converges toward 0 and for fixed $k \geq 2$, consider $\theta_k$ and $\xi^{(k)}$ be independent and respectively distributed as $\theta$ and $\xi$. Then we construct the sequence $(\theta_{k-1}, \xi^{(k-1)}), \ldots, (\theta_1, \xi^{(1)})$ as follows: for $j = k, k-1, \ldots, 2$,

$$(3.1) \qquad T_z^{(j)} \stackrel{(\text{def})}{=} \inf\{t : \xi_t^{(j)} \geq z\},$$

$$(3.2) \qquad \xi^{(j-1)} \stackrel{(\text{def})}{=} (\xi^{(j)}(T_{\log(x_{j-1}e^{-\theta_j}/x_j)}^{(j)} + t) - \xi^{(j)}(T_{\log(x_{j-1}e^{-\theta_j}/x_j)}^{(j)}), t \geq 0),$$

$$(3.3) \qquad \theta_{j-1} \stackrel{(\text{def})}{=} \xi^{(j)}(T_{\log(x_{j-1}e^{-\theta_j}/x_j)}^{(j)}) - \log(x_{j-1}e^{-\theta_j}/x_j).$$

So, this defines the law of $((\theta_1, \xi^{(1)}), \ldots, (\theta_k, \xi^{(k)}))$ on $(\mathbb{R}_+ \times \mathcal{D})^k$. We can check by induction and from (3.2), (3.3) and Proposition 1 that, for any $j = 1, \ldots, k$, $\theta_j$ and $\xi^{(j)}$ are independent and respectively distributed as $\theta$ and $\xi$.

Now let $k' > k$ and $\bar{\theta}_{k'}$ and $\bar{\xi}^{(k')}$ be independent and respectively distributed as $\theta$ and $\xi$. Then let us construct the sequence $((\bar{\theta}_1, \bar{\xi}^{(1)}), \ldots, (\bar{\theta}_{k'}, \bar{\xi}^{(k')}))$ on $(\mathbb{R}_+ \times \mathcal{D})^{k'}$ as above. Since $\bar{\theta}_k$ and $\bar{\xi}^{(k)}$ are independent and since $(\bar{\theta}_{k-1}, \bar{\xi}^{(k-1)}), \ldots, (\bar{\theta}_1, \bar{\xi}^{(1)})$ are constructed from $(\bar{\theta}_k, \bar{\xi}^{(k)})$ through the same way as the sequence $(\theta_{k-1}, \xi^{(k-1)}), \ldots, (\theta_1, \xi^{(1)})$ is constructed from $(\theta_k, \xi^{(k)})$, both sequences

$$((\bar{\theta}_1, \bar{\xi}^{(1)}), \ldots, (\bar{\theta}_k, \bar{\xi}^{(k)})) \quad \text{and} \quad ((\theta_1, \xi^{(1)}), \ldots, (\theta_k, \xi^{(k)}))$$

have the same law. For any $k \geq 1$, call $\mu^{(k)}$ the law of $((\theta_1, \xi^{(1)}), \ldots, (\theta_k, \xi^{(k)}))$ on the space $(\mathbb{R}_+ \times \mathcal{D})^k$. Then we proved that $(\mu^{(k)})$ is a projective family, that is, for any $k < k'$, the projection of the law $\mu^{(k')}$ on $(\mathbb{R}_+ \times \mathcal{D})^k$ corresponds to $\mu^{(k)}$. From Kolmogorov's theorem, there exists a unique probability law, say, $\mu$, on $(\mathbb{R}_+ \times \mathcal{D})^{\mathbb{N}}$ such that, for each $k$, $\mu^{(k)}$ is the projection of $\mu$ on $(\mathbb{R}_+ \times \mathcal{D})^k$. This defines the law of the sequence $(\theta_n, \xi^{(n)})$.

The relations (2.3) and (2.4) for $i = j$ are obvious, and for $i = j-1$, they correspond to (3.2) and (3.3). The case $i < j-1$ is easily obtained by induction. Let us check it for $i = j-2$. From (3.2),

$$\begin{aligned}
\xi^{(j-2)} = (\xi^{(j)}(T_{\log(x_{j-1}e^{-\theta_j}/x_j)}^{(j)} + T_{\log(x_{j-2}e^{-\theta_{j-1}}/x_{j-1})}^{(j-1)} + t) \\
- \xi^{(j)}(T_{\log(x_{j-1}e^{-\theta_j}/x_j)}^{(j)} + T_{\log(x_{j-2}e^{-\theta_{j-1}}/x_{j-1})}^{(j-1)}), t \geq 0),
\end{aligned}$$

and it is easy to see that

$$T_{\log(x_{j-1}e^{-\theta_j}/x_j)}^{(j)} + T_{\log(x_{j-2}e^{-\theta_{j-1}}/x_{j-1})}^{(j-1)} = T_{\log(x_{j-2}e^{-\theta_j}/x_j)}^{(j)}.$$

This gives (2.3), for $i = j-2$, and (2.4) follows.

Let $m \geq n \geq 1$. Since $\xi^{(m)}$ is independent of $\theta_m$, from (2.3) for $i = n$ and $j = m$, the process $\xi^{(n)}$ is independent of $\{\theta_m, (\xi_t^{(m)}, t \leq T_{\log(x_n e^{-\theta_m}/x_m)}^{(m)})\}$.

But from (2.4) for $i = m-1, m-2, \ldots, n$ and $j = m$, the variables $\theta_{m-1}, \ldots, \theta_n$ are functionals of $\{\theta_m, (\xi_t^{(m)}, t \leq T_{\log(x_n e^{-\theta_m}/x_m)}^{(m)})\}$. We deduce that $\xi^{(n)}$ is independent of $(\theta_k)_{k \geq n}$. Moreover, it follows directly from (2.4) that, for any $n$ and $m$ such that $n < m$,

$$
\begin{aligned}
(\theta_m, &\ldots, \theta_n, \theta_{n-1}, \ldots, \theta_1) \\
(3.4) \qquad &\overset{(a.s.)}{=} \left( \theta_m, \ldots, \theta_n, \xi_{T^{(n)}(\log x_{n-1}/x_n - \theta_n)}^{(n)} - \left( \log \frac{x_{n-1}}{x_n} - \theta_n \right), \ldots, \right. \\
&\left. \qquad \xi_{T^{(n)}(\log x_1/x_n - \theta_n)}^{(n)} - \left( \log \frac{x_1}{x_n} - \theta_n \right) \right).
\end{aligned}
$$

This shows that $(\theta_n)$ is a Markov chain, since in the right-hand side of (3.4), $\xi^{(n)}$, is independent of $(\theta_m, \theta_{m-1}, \ldots, \theta_n)$.  □

We shall see later that the tail sigma field $\mathcal{G} = \bigcap_n \sigma\{\theta_n, \theta_{n+1}, \ldots\}$ is trivial. Although, there should be more direct arguments, our proof bears on the construction of $X^{(0)}$, see Lemma 3. Note also that the Markov chain $(\theta_n)$ is homogeneous if and only if $x_{n-1}/x_n$ is constant, and in this case, it is stationary.

Let $\xi$ and the canonical process $X$ be related by (1.2) and define the first passage time process of $X$ by

$$
S_y = \inf\{t : X_t \geq y\}, \qquad y \geq 0.
$$

Observe that from our assumptions (see the beginning of Section 2.2), for any level $y \geq x$, $\mathbb{P}_x(S_y < \infty) = 1$.

LEMMA 1. *Let $x \leq y$ and set $z = \log(y/x)$. Then under $\mathbb{P}_x$, the process $(X_{S_y+t}, t \geq 0)$ admits the following (path by path) representation:*

$$
(3.5) \qquad (X(S_y + t), t \geq 0) = (y e^{\theta^{(z)}} \exp \bar{\xi}(\bar{\tau}(t e^{-\theta^{(z)}}/y)), t \geq 0),
$$

*where $\bar{\xi} \overset{(def)}{=} (\xi_{T_z+t} - \xi_{T_z}, t \geq 0)$ has the same law as $\xi$ and is independent of $\theta^{(z)} \overset{(def)}{=} \xi_{T_z} - z$ and $\bar{\tau}_t = \inf\{s : \int_0^s \exp \bar{\xi}_u \, du > t\}$.*

*In particular, when $\xi$ has no positive jumps, the overshoot $\theta^{(z)}$ is zero and $(X_{S_y+t}, t \geq 0)$ has law $\mathbb{P}_y$.*

PROOF. First observe that $S_y = x A_{T_z}$ and $X_{S_y} = x \exp \xi_{T_z}$, where $z = \log(y/x)$ and $T_z = \inf\{t : \xi_t \geq z\}$. Now, from (1.2), we have

$$
(X(S_y + t), t \geq 0) = (x \exp \xi_{\tau(A_{T_z} + t/x)}, t \geq 0).
$$

Then we can rewrite the time change as follows:

$$\tau(A_{T_z} + t/x) = \inf\{s \geq 0 : A_s \geq A_{T_z} + t/x\}$$
$$= \inf\left\{s \geq T_z : \int_{T_z}^s \exp(\xi_u - \xi_{T_z})\,du \geq (t/x)\exp(-\xi_{T_z})\right\}$$
$$= T_z + \bar{\tau}(t/X_{S_y}),$$

where $\bar{\tau}$ is the right continuous inverse of the exponential functional $\int_0^t \exp(\bar{\xi}_s)\,ds$, that is, $\bar{\tau}_t = \inf\{s : \int_0^s \exp(\bar{\xi}_u)\,du > t\}$, and $\bar{\xi} = (\xi_{T_z+t} - \xi_{T_z}, t \geq 0)$. Now we may write

$$(X(S_y + t), t \geq 0) = (x\exp\xi_{T_z+\bar{\tau}(t/X_{S_y})}, t \geq 0)$$
$$= (X_{S_y}\exp\bar{\xi}_{\bar{\tau}(t/X_{S_y})}, t \geq 0).$$

Finally, it follows from standard properties of Lévy processes that $\bar{\xi}$ has the same law as $\xi$ and is independent of $\xi_{T_z}$. □

COROLLARY 1. *Let $(\theta_n, \xi^{(n)})$, $n = 1, 2, \ldots$, be the sequence which is defined by Proposition 2 and recall the definition (2.5) of the sequence $X^{(\bar{x}_n)}$:*

$$X_t^{(\bar{x}_n)} = \bar{x}_n \exp\xi_{\tau^{(n)}(t/\bar{x}_n)}^{(n)}, \qquad t \geq 0, n \geq 1,$$

*with $\bar{x}_n = x_n e^{\theta_n}$ and $\tau_t^{(n)} \stackrel{(def)}{=} \inf\{s : \int_0^s \exp\xi_u^{(n)}\,du > t\}$. Recall also that $S^{(n-1)} = \inf\{t : X_t^{(\bar{x}_n)} \geq x_{n-1}\}$, then for every $n \geq 1$,*

$$(X^{(\bar{x}_n)}(S^{(n-1)} + t), t \geq 0) = X^{(\bar{x}_n)} \qquad a.s.$$

PROOF. Replacing $x$ by $\bar{x}_n = x_n e^{\theta_n}$ and $y$ by $x_{n-1}$ in the (path by path) identity of Lemma 1, we obtain [$\theta^{(z)}$ being defined in this lemma]

$$\theta^{(z)} = \xi^{(n)}(T_{\log(x_{n-1}/\bar{x}_n)}^{(n)}) - \log(x_{n-1}/\bar{x}_n) \qquad a.s.,$$

where $T_z^{(n)} = \inf\{t : \xi_t^{(n)} \geq z\}$. The right-hand side of this inequality is the r.v. $\theta_{n-1}$ which is defined in (2.4). We can also check that $\bar{\xi}$ defined in Lemma 1 is nothing but $\xi^{(n-1)}$ defined in (2.3) and the conclusion follows. □

The next result will be used in the sequel to check the fact that the construction (2.7) does not depend on the sequence $(x_n)$.

COROLLARY 2. *Let $(y_n)$ be another real decreasing sequence which tends to 0. For any $j$, let $n$ such that $x_n \leq y_j$ and let $V_j = \inf\{t : X_t^{(\bar{x}_n)} \geq y_j\}$, then we may write*

$$(3.6) \qquad (X^{(\bar{x}_n)}(V_j + t), t \geq 0) = (y_j e^{\tilde{\theta}_j} \exp\tilde{\xi}_{\tilde{\tau}^{(j)}(te^{-\tilde{\theta}_j}/y_j)}^{(j)}, t \geq 0),$$

*where*

$$\tilde{\xi}^{(j)} = (\xi^{(n)}(T^{(n)}_{\log(y_j/\bar{x}_n)} + t) - \xi^{(n)}(T^{(n)}_{\log(y_j/\bar{x}_n)}), t \geq 0,$$

$$\tilde{\theta}_j = \xi^{(n)}(T^{(n)}_{\log(y_j/\bar{x}_n)}) - \log(y_j/\bar{x}_n)$$

*and* $\tilde{\tau}^{(j)}_t = \inf\{s : \int_0^s \exp \tilde{\xi}^{(j)}_u \, du > t\}$.

*The sequence* $(\tilde{\theta}_j, \tilde{\xi}^{(j)})$, $j \geq 1$, *may be defined the same way as in Proposition 2 with respect to the sequence* $(y_n)$. *That is,* $\tilde{\theta}_j$ *and* $\tilde{\xi}^{(j)}$ *are independent and respectively distributed as* $\theta$ *and* $\xi$, *and for* $1 \leq i \leq j$,

$$(3.7) \qquad \tilde{\xi}^{(i)} = (\tilde{\xi}^{(j)}(\tilde{T}^{(j)}_{\log(y_i e^{-\tilde{\theta}_j}/y_j)} + t) - \tilde{\xi}^{(j)}(\tilde{T}^{(j)}_{\log(y_i e^{-\tilde{\theta}_j}/y_j)}), t \geq 0,$$

$$(3.8) \qquad \tilde{\theta}_i = \tilde{\xi}^{(j)}(\tilde{T}^{(j)}_{\log(y_i e^{-\tilde{\theta}_j}/y_j)}) - \log(y_i e^{-\tilde{\theta}_j}/y_j),$$

*where* $\tilde{T}^{(j)}_z = \inf\{t : \tilde{\xi}^{(j)}_t \geq z\}$. *In particular, the law of* $(\tilde{\theta}_j, \tilde{\xi}^{(j)})$, $j \geq 1$, *does not depend on the sequence* $(x_n)$.

PROOF. The first part concerning the definitions of $\tilde{\xi}^{(j)}$, $\tilde{\theta}_j$ and $\tilde{\tau}^{(j)}$ is a straightforward consequence of Lemma 1.

Now we prove that the sequence $(\tilde{\theta}_j, \tilde{\xi}^{(j)})$ may be defined as in Proposition 2. First, it is clear from the independence between $\theta_n$ and $\xi^{(n)}$ and from the stationarity of the process $(\xi_{T_{z-\theta}} - (z - \theta), z \geq 0)$ that $\tilde{\xi}^{(j)}$ and $\tilde{\theta}^{(j)}$ are independent and respectively distributed as $\xi$ and $\theta$. It remains to check that $(\tilde{\theta}_j, \tilde{\xi}^{(j)})$ verifies (3.7) and (3.8). Let $i \leq j$ so that $y_j \leq y_i$, then, from the statement, $\tilde{\xi}^{(j)}$ and $\tilde{\xi}^{(i)}$ are respectively defined by $\tilde{\xi}^{(j)} = (\xi^{(n)}_{T^{(n)}(\log(y_j/\bar{x}_n))+t} - \xi^{(n)}_{T^{(n)}(\log(y_j/\bar{x}_n))}, t \geq 0)$ and $\tilde{\xi}^{(i)} = (\xi^{(n')}_{T^{(n')}(\log(y_i/\bar{x}_{n'}))+t} - \xi^{(n')}_{T^{(n')}(\log(y_i/\bar{x}_{n'}))}, t \geq 0)$ for some indices $n$ and $n'$ such that $x_n \leq y_j$ and $x_{n'} \leq y_i$. But we can check from Corollary 1 that since $x_n \leq y_i$ and $x'_n \leq y_i$, we have $X^{(\bar{x}_n)}_{V_i+} = X^{(\bar{x}_{n'})}_{V'_i+}$, where $V'_j = \inf\{t : X^{(\bar{x}_{n'})}_t \geq y_j\}$, so that $\tilde{\xi}^{(i)}$ may be defined with respect to $X^{(\bar{x}_n)}$, that is,

$$(X^{(\bar{x}_n)}(V_i + t), t \geq 0) = (y_i e^{\tilde{\theta}_i} \exp \tilde{\xi}^{(i)}_{\tilde{\tau}^{(i)}(te^{-\tilde{\theta}_i}/y_i)}, t \geq 0),$$

*and, in particular, we have*

$$\tilde{\xi}^{(i)} = (\xi^{(n)}(T^{(n)}_{\log(y_i/\bar{x}_n)} + t) - \xi^{(n)}(T^{(n)}_{\log(y_i/\bar{x}_n)}), t \geq 0).$$

On the other hand, one easily checks that

$$T^{(n)}_{\log(y_i/\bar{x}_n)} = T^{(n)}_{\log(y_j/\bar{x}_n)} + \tilde{T}^{(j)}_{\log(y_i e^{-\tilde{\theta}_j}/y_j)},$$

hence,

$$
\tilde{\xi}^{(i)} = (\xi^{(n)}(T^{(n)}_{\log(y_j/\bar{x}_n)} + \tilde{T}^{(j)}_{\log(y_i e^{-\tilde{\theta}_j}/y_j)} + t)
$$

$$
- \xi^{(n)}(T^{(n)}_{\log(y_j/\bar{x}_n)} + \tilde{T}^{(j)}_{\log(y_i e^{-\tilde{\theta}_j}/y_j)}), t \geq 0)
$$

$$
= (\tilde{\xi}^{(j)}(\tilde{T}^{(j)}_{\log(y_i e^{-\tilde{\theta}_j}/y_j)} + t) - \tilde{\xi}^{(j)}(\tilde{T}^{(j)}_{\log(y_i e^{-\tilde{\theta}_j}/y_j)}), t \geq 0),
$$

which is identity (3.7), and identity (3.8) follows. $\square$

PROOF OF THEOREM 1.   Since $X^{(0)}$ is obtained as the concatenation of the processes $(X^{(\bar{x}_n)}_{\Sigma_n - t}, \Sigma_n \leq t \leq \Sigma_{n-1})$, we need conditions which insure that the sum $\Sigma_n = \sum_{k=n}^{\infty} S^{(k)}$ is a.s. finite. Moreover, for $X^{(0)}$ to be issued from 0 in a "continuous" way, that is, otherwise than by a jump, it is also necessary to have $\lim_n \bar{x}_n = 0$. These are the objects of the next lemma.

LEMMA 2.   *Assume that* (H) *holds:*

 (i)  $\lim_{n \to +\infty} \bar{x}_n = 0$, *a.s.*
 (ii) *For any* $n \geq 1$, $\Sigma_n > 0$, *a.s.*
 (iii) $\Sigma_n < +\infty$, *a.s. if and only if* (2.11) *holds.*

PROOF.   To prove part (i), it suffices to observe that $\bar{x}_n = x_n e^{\theta_n}$ is a nonnegative decreasing sequence [indeed, $\bar{x}_{n-1} = X^{(\bar{x}_n)}(S^{(n-1)})$], so it converges almost surely. Moreover, since the $\theta_n$'s are identically distributed, the limit law of the sequence $(x_n e^{\theta_n})$ is the Dirac mass at 0.

To prove (ii), first note that $\Sigma_n = 0$ a.s. if and only if, for each $k \geq n$, $S^{(k)} = 0$, which is also equivalent to $X^{(\bar{x}_k)}_0 = x_k e^{\theta_k} \geq x_{k-1}$. Indeed, $X^{(\bar{x}_k)}$ is right continuous, so it has no jump at 0. Hence, we have

$$
P(\Sigma_n = 0) = P(x_j e^{\theta_j} \geq x_{j-1}, \text{for all } j \geq n)
$$

$$
= \lim_{k \to +\infty} P(x_j e^{\theta_j} \geq x_{j-1}, \text{for all } j = n, \ldots, n+k).
$$

But from (3.4), $\{x_j e^{\theta_j} \geq x_{j-1}, \text{for all } j = n, \ldots, n+k\} = \{x_{n+k} e^{\theta_{n+k}} \geq x_{n-1}\}$ so that

$$
P(\Sigma_n = 0) = \lim_{k \to +\infty} P(x_{n+k} e^{\theta_{n+k}} \geq x_{n-1}),
$$

and this limit is 0 since the $\theta_i$'s have the same nondegenerate law and $(x_n)$ tends to 0.

Now we prove (iii). We start by proving that the convergence of $\Sigma_n = \sum_{k \geq n} S^{(k)}$ does not depend on the choice of the sequence $(x_n)$. Let

$$
\tilde{X}^{(\bar{y}_n)} = (y_n e^{\tilde{\theta}_n} \exp \tilde{\xi}^{(n)}_{\tilde{\tau}^{(n)}(te^{-\tilde{\theta}_n}/y_n)}, t \geq 0),
$$

with $\bar{y}_n = y_n e^{\tilde{\theta}_n}$, be the processes which are defined in (3.6) of Corollary 2, with respect to a sequence $(y_n)$ which decreases toward 0. Set $\widetilde{S}^{(n-1)} = \inf\{t : \widetilde{X}_t^{(\bar{y}_n)} \geq y_{n-1}\}$, then we see from the identity (3.6) that if $x_n \leq y_j \leq y_{j-1} \leq x_k$, then $S^{(n-1)} + \cdots + S^{(k)} \geq \widetilde{S}^{(j-1)}$. More generally, for any $i, j$, there exist $n, m$ such that

$$(3.9) \qquad \sum_n^m S^{(k)} \geq \sum_i^j \widetilde{S}^{(k)}.$$

Conversely, one can check exactly as in Corollary 2 that, for $y_n \leq x_j$, if $\widetilde{V}_j = \inf\{t : \widetilde{X}_t^{(\bar{y}_n)} \geq x_j\}$, then we have

$$(\widetilde{X}^{(\bar{y}_n)}(\widetilde{V}_j + t), t \geq 0) = (x_j e^{\theta_j} \exp \xi_{\tau^{(j)}(te^{-\theta_j}/x_j)}^{(j)}, t \geq 0),$$

which shows that if $y_n \leq x_j \leq x_{j-1} \leq y_k$, then $\widetilde{S}^{(n-1)} + \cdots + \widetilde{S}^{(k)} \geq S^{(j-1)}$. Then as above, for any $i, j$, there exist $n, m$ such that

$$(3.10) \qquad \sum_n^m \widetilde{S}^{(k)} \geq \sum_i^j S^{(k)}.$$

Inequalities (3.9) and (3.10) prove that $\sum_{k \geq n} \widetilde{S}^{(k)}$ is finite if and only if $\sum_{k \geq n} S^{(k)}$ is finite. Hence, the convergence of the sum $\Sigma_n = \sum_{k \geq n} S^{(k)}$ does not depend on the sequence $(x_n)$ and we can consider any particular sequence in the sequel of this proof.

Now observe that $S^{(n-1)} = x_n e^{\theta_n} \int_0^{T_{\nu_n}^{(n)}} \exp \xi_s^{(n)} \, ds$, where $\nu_n = \log(x_{n-1} e^{-\theta_n}/x_n)$ and $T_z^{(n)} = \inf\{t \geq 0 : \xi_t^{(n)} \geq z\}$. Note also that $T_{\nu_n}^{(n)} = 0$ whenever $x_n e^{\theta_n} \geq x_{n-1}$. Hence,

$$x_n \int_0^{T_{\nu_n}^{(n)}} \exp \xi_s^{(n)} \, ds \leq S^{(n-1)} \leq x_{n-1} \int_0^{T_{\nu_n}^{(n)}} \exp \xi_s^{(n)} \, ds \qquad \text{a.s.}$$

These inequalities imply that, for any sequence $(z_n)$ of positive reals,

$$(3.11) \quad \begin{aligned} & x_n \mathbb{1}_{\{x_n e^{\theta_n} < z_n\}} \int_0^{T^{(n)}(\log x_{n-1}/z_n)} \exp \xi_s^{(n)} \, ds \\ & \qquad \leq S^{(n-1)} \leq x_{n-1} \int_0^{T^{(n)}(\log x_{n-1}/x_n)} \exp \xi_s^{(n)} \, ds \qquad \text{a.s.} \end{aligned}$$

Note also that, for any $r > 1$,

$$(3.12) \quad E\Big(\log^+ \int_0^{T_1} \exp \xi_s \, ds\Big) < \infty \quad \Longleftrightarrow \quad \sum P\Big(\int_0^{T_1} \exp \xi_s \, ds > r^n\Big) < \infty.$$

By taking $r = e^{n/2}$ in (3.12), and from Borel–Cantelli's lemma, we obtain that if $E(\log^+ \int_0^{T_1} \exp \xi_s \, ds) < \infty$, then $e^{-n} \int_0^{T_1^{(n)}} \exp \xi_s^{(n)} \, ds < e^{-n/2}$, a.s. for $n$ sufficiently large. So by choosing $x_n = e^{-n}$ in (3.11), we obtain that $\sum S^{(n)} < \infty$, a.s.

Conversely, suppose that $E(\log^+ \int_0^{T_1} \exp \xi_s \, ds) = +\infty$. Take $r > 1$, $x_n = r^{-n}$ and $z_n = (1 - r^{-1/2})r^{-n} + r^{-n+1/2}$ in (3.11) and set

$$A_n = \left\{ \int_0^{T^{(n)}(\log x_{n-1}/z_n)} \exp \xi_s^{(n)} \, ds > r^n \right\},$$

$$B_n = \left\{ x_n \mathbb{1}_{\{x_n e^{\theta_n} < z_n\}} \int_0^{T^{(n)}(\log x_{n-1}/z_n)} \exp \xi_s^{(n)} \, ds > 1 \right\} = A_n \cap \{x_n e^{\theta_n} < z_n\}.$$

Let $r > 1$ be sufficiently large to have $\frac{x_{n-1}}{z_n} = \frac{r}{(1 - r^{-1/2}) + r^{1/2}} > e$, then, from (3.12), $\sum P(\int_0^{T^{(n)}(\log x_{n-1}/z_n)} \exp \xi_s^{(n)} \, ds > r^n) = \infty$, so that from the independence between $\theta_n$ and $\xi^{(n)}$, $\sum P(B_n) = \sum P(\theta < \log(1 - r^{-1/2} + r^{1/2}))P(A_n) = +\infty$, since, from (2.2), $P(\theta < \eta) > 0$, for all $\eta > 0$. Moreover, for all $n, m \geq 1$ such that $n \neq m$,

$$P(B_n \cap B_m) \leq P(A_n \cap A_m) = \frac{P(B_n)P(B_m)}{P(\theta < \log(1 - r^{-1/2} + r^{1/2}))^2}.$$

Hence, from a well-known extension of Borel–Cantelli's lemma (see, e.g., [14], page 317), we have $P(\limsup_n B_n) > P(\theta < \log(1 - r^{-1/2} + r^{1/2}))^2$. Then for all $\varepsilon > 0$, there exists $r_0$ such that, for any $r > r_0$, $P(\limsup_n B_n) > 1 - \varepsilon$ (note that $B_n$ depends on $r$). But $\limsup_n B_n \subset \limsup_n A_n \subset \{\sum_n S^{(n)} = +\infty\}$ and the probability of the last event does not depend on $r$. Therefore, $P(\sum_n S^{(n)} = +\infty) = 1$. $\quad\square$

Lemma 2 shows that the definition of $X^{(0)}$, that is,

$$X_t^{(0)} = \begin{cases} X_{t-\Sigma_1}^{(\bar{x}_1)}, & t \in [\Sigma_1, \infty), \\ X_{t-\Sigma_2}^{(\bar{x}_2)}, & t \in [\Sigma_2, \Sigma_1), \\ \vdots & \\ X_{t-\Sigma_n}^{(\bar{x}_n)}, & t \in [\Sigma_n, \Sigma_{n-1}), \\ \vdots & \end{cases}$$

makes sense. Indeed, since $\lim_n \Sigma_n = 0$, the process $X^{(0)}$ is well defined on $(0, \infty)$. Moreover, since $\Sigma_n > 0$, a.s. for any $n$ and $\lim_n x_n e^{\theta_n} = 0$, a.s., we have $\lim_{t \to 0} X_t^{(0)} = 0$. Now we put $X_0^{(0)} = 0$ so that $X^{(0)}$ is a càdlàg process on $[0, \infty)$, which is positive on $(0, \infty)$.

As in the proof of the previous lemma, let

$$\widetilde{X}^{(\bar{y}_n)} = (y_n e^{\tilde{\theta}_n} \exp \xi_{\tilde{\tau}^{(n)}(te^{-\tilde{\theta}_n}/y_n)}^{(n)}, t \geq 0)$$

be the processes which are defined in Corollary 2 and recall $\widetilde{S}^{(n-1)} = \inf\{t : \widetilde{X}^{(\bar{y}_n)} \geq y_{n-1}\}$ and $\widetilde{\Sigma}_n = \sum_{k \geq n} \widetilde{S}^{(k)}$. Let $\bar{y}_n = y_n e^{\hat{\theta}_n}$, then we see from the construction (2.7) and from (3.6) in Corollary 2 that $X^{(0)}$ may be represented as the concatenation of the processes $(\widetilde{X}^{(\bar{y}_n)}_{\widetilde{\Sigma}_n - t}, \widetilde{\Sigma}_n \leq t \leq \widetilde{\Sigma}_{n-1})$ which are defined as $(X^{(\bar{x}_n)}_{\Sigma_n - t}, \Sigma_n \leq t \leq \Sigma_{n-1})$, in the same way as in (2.7), but with respect to the sequence $(y_n)$. The law of these processes does not depend on the sequence $(x_n)$, hence, neither does the law of $X^{(0)}$. This proves (i) of Theorem 1.

Part (ii) follows from the path properties of the processes $X^{(\bar{x}_n)}$ and is straightforward.

We see from the construction (2.7) of the process $X^{(0)}$ that the process $(kX^{(0)}_{k^{-1}t}, t \geq 0)$ is obtained from the concatenation of the processes $(\widetilde{X}^{(\bar{y}_n)}_{\widetilde{\Sigma}_n - t}, \widetilde{\Sigma}_n \leq t \leq \widetilde{\Sigma}_{n-1})$, with the particular sequence $y_n = kx_n$ [note that, in this particular case, $(\tilde{\theta}_n, \tilde{\xi}^{(n)}) = (\theta_n, \xi^{(n)})$]. Hence, from (i) of this theorem, the process $(kX^{(0)}_{k^{-1}t}, t \geq 0)$ has the same law as $X^{(0)}$ and (iii) is proved.

It remains to prove (iv) of Theorem 1, that is, the Markov property of $X^{(0)}$. To do so, we need the following lemma.

LEMMA 3. *The process $X^{(0)}$ is independent of the tail-sigma field of $(\theta_n)$, that is, $\mathcal{G} \stackrel{(def)}{=} \bigcap_n \sigma\{\theta_i, i \geq n\}$. Consequently, $\mathcal{G}$ is trivial.*

PROOF. As we already showed, it follows from Corollary 1 and construction (2.7) that $X^{(\bar{x}_n)} = (X^{(0)}_{\Sigma_n + t}, t \geq 0)$. We derive from this identity that $\lim_{n \to \infty} X^{(\bar{x}_n)} = X^{(0)}$, a.s. on the space $\mathcal{D}$ of càdlàg trajectories. Let $\mathcal{G}_n = \sigma\{\theta_n, \theta_{n+1}, \ldots\}$ and $\mathcal{G} = \bigcap_n \mathcal{G}_n$. By construction, $X^{(\bar{x}_n)}$ is a functional of $\xi^{(n)}$ and $\theta_n$, then, from Proposition 1, it is clear that the law of this process conditionally on $(\theta_n, \theta_{n+1}, \ldots)$ is the same as its law conditionally on $\theta_n$, so that, from the Markov property, for any bounded and continuous functional $H$ on $\mathcal{D}$, we have

$$(3.13) \quad E(H(X^{(\bar{x}_n)})|\mathcal{G}_n) = E(H(X^{(\bar{x}_n)})|\theta_n) = \mathbb{E}_{x_n e^{\theta_n}}(H) \xrightarrow[n \to \infty]{\text{(a.s.)}} E(H(X^{(0)})|\mathcal{G}).$$

Now fix $k \geq 1$ and let $f_1, \ldots, f_k$ be $k$ bounded measurable functions, then it follows from (3.13) that

$$(3.14) \quad E(f_1(\theta_1) \cdots f_k(\theta_k) E(H(X^{(0)})|\mathcal{G})) = \lim_n E(f_1(\theta_1) \cdots f_k(\theta_k) \mathbb{E}_{x_n e^{\theta_n}}(H)).$$

From our construction of $\theta_1, \ldots, \theta_k$ and, more precisely, from (3.4), we have

$$E(f_1(\theta_1) \cdots f_k(\theta_k) \mathbb{E}_{x_n e^{\theta_n}}(H))$$
$$= E(f_1[\xi^{(n)}_{T^{(n)}(\log x_1 e^{-\theta_n}/x_n)} - (\log x_1 e^{-\theta_n}/x_n)]$$

$$\times \cdots \times f_k[\xi^{(n)}_{T^{(n)}(\log x_k e^{-\theta_n}/x_n)} - (\log x_k e^{-\theta_n}/x_n)]\mathbb{E}_{x_n e^{\theta_n}}(H)),$$

where we recall that $\xi^{(n)}$ and $\theta_n$ are independent. Denote by $E_z$ the law of the process $\xi^{(n)} + z$, $z \in \mathbb{R}$. Then since the $\theta_n$'s have the same law, the second member of the above equality may be written as

$$E(E_\theta(f_1[\xi^{(n)}_{T^{(n)}(\log x_1/x_n)} - (\log x_1/x_n)]$$

$$\times \cdots \times f_k[\xi^{(n)}_{T^{(n)}(\log x_k/x_n)} - (\log x_k/x_n)])\mathbb{E}_{x_n e^{\theta}}(H)),$$

and from Proposition 1,

$$(3.15) \quad \begin{aligned} E_\theta(f_1[\xi^{(n)}_{T^{(n)}(\log x_1/x_n)} - (\log x_1/x_n)] \cdots f_k[\xi^{(n)}_{T^{(n)}(\log x_k/x_n)} - (\log x_k/x_n)]) \\ \overset{\text{(a.s.)}}{\underset{n\to\infty}{\longrightarrow}} E(f_1(\theta_1) \cdots f_k(\theta_k)). \end{aligned}$$

We then deduce from (3.13), (3.14), (3.15) and Lebesgue's theorem of dominated convergence that

$$E(f_1(\theta_1) \cdots f_k(\theta_k)E(H(X^{(0)})|\mathcal{G})) = E(f_1(\theta_1) \cdots f_k(\theta_k))E(H(X^{(0)})).$$

As a consequence, we have $E(E(H(X^{(0)})|\mathcal{G})|\sigma(\theta_1, \ldots, \theta_k)) = E(H(X^{(0)}))$, for any $k$, so that $E(E(H(X^{(0)})|\mathcal{G})|\sigma(\theta_1, \theta_2 \ldots)) = E(H(X^{(0)}))$, hence, $E(H(X^{(0)})|\mathcal{G}) = E(H(X^{(0)}))$ since $\mathcal{G} \subset \sigma(\theta_1, \theta_2 \ldots)$. It proves that $X^{(0)}$ is independent of $\mathcal{G}$. On the other hand, since $\Sigma_n = \inf\{t : X_t^{(0)} \geq x_n\}$ and $X^{(0)}_{\Sigma_n} = X^{(\bar{x}_n)}_0 = x_n e^{\theta_n}$, then the sequence $(\theta_n)$ may be generated by the process $X^{(0)}$. In particular, $\mathcal{G}$ is a sub-$\sigma$ field of the $\sigma$-field generated by $X^{(0)}$, hence, it is trivial. □

To prove that $X^{(0)}$ is a strong Markov process, note that, from (3.13), and Lemma 3, we obtain

$$(3.16) \quad \mathbb{E}_{\bar{x}_n}(H) \overset{\text{(a.s.)}}{\underset{n\to\infty}{\longrightarrow}} E(H(X^{(0)})),$$

for any bounded, continuous functional $H$ defined on $\mathcal{D}$. So the Markov property of the process $X^{(0)}$ is inherited from the same property for the family $\mathbb{P}_x$, $x > 0$. This proves part (iv) and ends the proof of Theorem 1. □

PROOF OF PROPOSITION 3. Suppose first that $E(\xi_1) > 0$. It has been shown in [2] and [4] that, for fixed $t \geq 0$, the law of $X_t$ under $\mathbb{P}_x$ converges weakly as $x \to 0$ toward a nondegenerate probability law. It is clear from (3.16) that this entrance law is the law of $X_t^{(0)}$, so the result is due to [2] and [4] in that case.

Now suppose that $E(\xi_1) = 0$ and consider the process $X^{(0,\lambda)}$ which is defined as in (2.7), but relatively to the Lévy process $\xi^{(\lambda)} \stackrel{\text{(def)}}{=} (\xi_t + \lambda t, t \geq 0)$ and its corresponding limit overshoot $\theta^{(\lambda)}$. [Note that, from Section 2.1, conditions of Theorem 1 are satisfied for $(\xi_t + \lambda t, t \geq 0)$ so that such a construction is possible.] Comparing (2.9) and (2.10), it appears that, to prove Proposition 3, we only have to check that, for all $\mu \geq 0$ and for all $t \geq 0$,

$$\lim_{\lambda \to 0} E(\exp(-\mu X_t^{(0,\lambda)})) = E(\exp(-\mu X_t^{(0)})).$$

Now, let $(x_n)$ be any sequence which is as in Proposition 2 and let $(\theta_n^{(\lambda)}, \xi^{(n,\lambda)})$,

$$X_t^{(n,\lambda)} = x_n e^{\theta_n^{(\lambda)}} \exp \xi_{\tau^{(n,\lambda)}(te^{-\theta_n^{(\lambda)}}/x_n)}^{(n,\lambda)}, \qquad t \geq 0, n \geq 1,$$

$\tau^{(n,\lambda)}$, and $S^{(n,\lambda)}$ be defined as in (2.3), (2.4), (2.5) and (2.6), but relative to the laws of $\xi^{(\lambda)}$ and its corresponding limit overshoot $\theta^{(\lambda)}$. Put $\Sigma_n^{(\lambda)} = \sum_{k \geq n} S^{(k,\lambda)}$. Then it follows from construction (2.7), stated in Theorem 1, and the Markov property that

$$
\begin{aligned}
(3.17) \quad E(\exp(-\mu X_t^{(0,\lambda)})) &= \sum_{n \geq 1} E(\exp(-\mu X_{t-\Sigma_n^{(\lambda)}}^{(n,\lambda)}) \mathbb{1}_{\{t \in [\Sigma_n^{(\lambda)}, \Sigma_{n-1}^{(\lambda)})\}}) \\
&= \sum_{n \geq 1} E(\mathbb{P}_{x_n e^{\theta_n^{(\lambda)}}}^{(\lambda)}(\exp(-\mu X_t) \mathbb{1}_{\{t < S_{x_{n-1}}\}})),
\end{aligned}
$$

where $\mathbb{P}_x^{(\lambda)}$ is the law of the process defined in (1.2), with respect to $\xi^{(\lambda)}$, $X$ is the canonical process on $\mathcal{D}$ and $S_{x_{n-1}}$ is its first passage time above the state $x_{n-1}$. On one hand, a result due to Doney [9] asserts that

$$(3.18) \qquad \theta^{(\lambda)} \xrightarrow{\text{(w)}} \theta \qquad \text{as } \lambda \text{ tends to } 0.$$

On the other hand, it is not difficult to check from (1.2) that the family $(\mathbb{P}_x^{(\lambda)})$ is weakly continuous in $(\lambda, x)$, on the space $[0, \infty) \times (0, \infty)$ and that the set of discontinuities of the functional $S_{x_{n-1}}$ is negligible under any of the probability measures $(\mathbb{P}_x^{(\lambda)})$, $(\lambda, x) \in [0, \infty) \times (0, \infty)$. Indeed, this set is

$$\{\omega \in \mathcal{D} : \exists \varepsilon > 0, X(S_{x_{n-1}} + t)(\omega) = x_{n-1}, \forall t \in [0, \varepsilon]\}$$

$$\cup \{\omega \in \mathcal{D} : \exists \varepsilon > 0, X(S_{x_{n-1}})(\omega) = x_{n-1},$$

$$X(S_{x_{n-1}} + t)(\omega) < x_{n-1}, \forall t \in (0, \varepsilon]\},$$

but the latter is negligible from (1.2) and the path properties of (nonlattice) Lévy processes. Hence, $\mathbb{E}_x^{(\lambda)}(\exp(-\mu X_t) \mathbb{1}_{\{t < S_{x_{n-1}}\}})$ is bounded and continuous in $(\lambda, x)$, on the space $[0, \infty) \times (0, \infty)$ and from (3.18), for any $n \geq 1$, since $x_n > 0$,

$$E(\mathbb{E}_{x_n e^{\theta_n^{(\lambda)}}}^{(\lambda)}(\exp(-\mu X_t) \mathbb{1}_{\{t < S_{x_{n-1}}\}})) \longrightarrow E(\mathbb{E}_{x_n e^{\theta_n}}(\exp(-\mu X_t) \mathbb{1}_{\{t < S_{x_{n-1}}\}})),$$

as $\lambda$ tends to 0, so that from (3.17) and Fatou's lemma we have, for all $\mu \geq 0$,

$$\lim_{\lambda \to 0} E(\exp(-\mu X_t^{(0,\lambda)})) = \sum_{n \geq 1} E(\mathbb{E}_{x_n e^{\theta_n}}(\exp(-\mu X_t)\mathbb{1}_{\{t < S_{x_{n-1}}\}}))$$

$$= E(\exp(-\mu X_t^{(0)})),$$

and the conclusion follows.  □

PROOF OF THEOREM 2. We first suppose that (H) holds and that $E(\log \int_0^{T_1} \exp \xi_s \, ds) < +\infty$. Let $H$ be any bounded, continuous functional $H$ which is defined on $\mathcal{D}$. Since the $\theta_n$'s have the same law, we deduce from (3.16) that

$$\mathbb{E}_{x_n e^{\theta}}(H) \xrightarrow[n \to \infty]{(\mathrm{P})} E(H(X^{(0)})).$$

It proves that $\mathbb{E}_{x_{i_n} e^{\theta}}(H) \xrightarrow[n \to \infty]{(\mathrm{a.s.})} E(H(X^{(0)}))$, for a subsequence $(x_{i_n})$. Let $\omega \in \Omega$ be such that $k \stackrel{(\mathrm{def})}{=} e^{\theta(\omega)}$ satisfies $\mathbb{E}_{k x_{i_n}}(H) \longrightarrow E(H(X^{(0)}))$, then, from the scaling property under $\mathbb{P}_x$, $x > 0$ and under $\mathbb{P}_0$ stated, respectively, in (1.1) and (2.8), we can replace $X^{(0)}$ by $(k X_{k^{-1}t}^{(0)}, t \geq 0)$ in the above convergence in order to obtain

$$\mathbb{E}_{x_{i_n}}(H) \longrightarrow E(H(X^{(0)})) \qquad \text{as } n \to +\infty.$$

We proved that, for any decreasing sequence $(x_n)$ which converges toward 0, there exists a subsequence $(x_{i_n})$ such that the above convergence holds, hence, it holds for any decreasing sequence which converges toward 0.

Conversely, suppose that the family of probability measures $(\mathbb{P}_x)$ converges weakly as $x \to 0$ toward a nondegenerate probability measure $\mathbb{P}_0$ on $\mathcal{D}$. Necessarily, $\mathbb{P}_0$ is the law of a pssMp with index 1, which starts from 0 and never comes back to 0 and whose lim sup is infinite. We denote by $X^{(0)}$ a process on $\mathcal{D}$ whose law is $\mathbb{P}_0$. [Here we no longer suppose the validity of the construction (2.7).] Let $(x_n)$ be a real decreasing sequence which tends to 0 and define $\Sigma_n = \inf\{t : X_t^{(0)} \geq x_n\}$. We see from the Markov property and (1.2) that, for any $n \geq 1$,

$$(3.19) \qquad (X^{(0)}(\Sigma_n + t), t \geq 0) \stackrel{(\mathrm{d})}{=} (X_{\Sigma_n}^{(0)} \exp \xi_{\tau^{(n)}(t/X_{\Sigma_n}^{(0)})}^{(n)}, t \geq 0),$$

where $\xi^{(n)}$ is distributed as $\xi$ and is independent of $X_{\Sigma_n}^{(0)}$, and $\tau^{(n)}$ is as usual. Let $n \geq 1$; from Lemma 1, the value of the process $(X_{\Sigma_n}^{(0)} \exp \xi_{\tau^{(n)}(t/X_{\Sigma_n}^{(0)})}^{(n)}, t \geq 0)$ at its first passage time above $x_1$ is

$$X_{\Sigma_n}^{(0)} \exp \xi^{(n)}(T_{\log(x_1/X_{\Sigma_n}^{(0)})}) = x_1 \exp(\xi^{(n)}(T_{\log(x_1/X_{\Sigma_n}^{(0)})}) - \log(x_1/X_{\Sigma_n}^{(0)})).$$

So from (3.19), we deduce that

$$(3.20) \qquad X^{(0)}(\Sigma_1) \overset{(d)}{=} x_1 \exp(\xi^{(n)}(T_{\log(x_1/X^{(0)}_{\Sigma_n})}) - \log(x_1/X^{(0)}_{\Sigma_n})).$$

Indeed, since $x_1 \geq x_n$, the value of the process $X^{(0)}$ at its first passage time above $x_1$ is the same as the value of the process $(X^{(0)}(\Sigma_n + t), t \geq 0)$ at its first passage time above $x_1$. On the other hand, the scaling property of the process $X^{(0)}$ implies that the law of the r.v.'s $x_n^{-1} X^{(0)}_{\Sigma_n}$, $n \geq 1$, is the same and does not depend on the sequence $(x_n)$. Moreover, these r.v.'s are a.s. positive. Let $\theta$ be an r.v. which is distributed as $\log(x_n^{-1} X^{(0)}_{\Sigma_n})$ and independent of the sequence $(\xi^{(n)})$, then we deduce from the fact that the law of $x_n^{-1} X^{(0)}_{\Sigma_n}$ does not depend on $n$ and from (3.20) that the law of $\xi^{(n)}(T_{\log(x_1/x_n)-\theta}) - (\log(x_1/x_n) - \theta)$ does not depend on $n$ and is the same as this of $\theta$. The same result is true for any sequence $(x_n)$ which decreases toward 0, hence, we have proved that the process $(\xi^{(n)}(T_{z-\theta}) - (z - \theta), z \geq 0)$ is stationary and Proposition 1 allows us to conclude that $\xi^{(n)}(T_z) - z$ converges weakly as $z$ tends to $+\infty$ toward the law of $\theta$.

Finally, recall that, from (3.19), we have

$$(X^{(0)}(\Sigma_n + t), t \geq 0) \overset{(d)}{=} (x_n e^{\theta_n} \exp \xi^{(n)}_{\tau^{(n)}(t/x_n e^{\theta_n})}, t \geq 0),$$

where $\theta_n = \log(x_n^{-1} X^{(0)}_{\Sigma_n})$ has the same law as $\theta$. Therefore, $X^{(0)}$ may be constructed as in (2.7) and the first passage time $\Sigma_n$ admits the same decomposition as in Theorem 1. This r.v. is obviously finite and Lemma 1(iii) shows that $E(\log^+ \int_0^{T_1} \exp \xi_s^{(n)} ds) < \infty$.  $\square$

REMARK.   When (H) holds but $E(\log^+ \int_0^{T_1} \exp \xi_s ds) = +\infty$, we see from Lemma 2 that $\Sigma_n = +\infty$. Following the previous proofs, it shows that the first passage time of the process $(X, \mathbb{P}_x)$ over $y$ (i.e., $S_y = x A_{T_z}$ in the notation of Lemma 1) tends almost surely toward $+\infty$, as $x$ tends to 0. It means that the process $(X, \mathbb{P}_x)$ converges almost surely toward the process $X^{(0)} \equiv 0$.

## 4. Examples.
The aim of this section is to present some examples of pssMp for which the conditions of Theorems 1 and 2 are satisfied.

Recall that the only positive self-similar Markov processes which are continuous are Bessel processes raised to any nonnegative power and multiplied by any constant. Indeed, since the only continuous Lévy processes are Brownian motions with drift multiplied by constant, this observation is a direct consequence of Lamperti representation (1.2). Recall also that the Bessel process of dimension $\mu \geq 0$ and starting from $x \geq 0$ is the diffusion $R$ whose square satisfies the stochastic differential equation

$$(4.1) \qquad R_t^2 = x^2 + 2 \int_0^t R_s \, dB_s + \mu t, \qquad t \geq 0,$$

where $B$ is the standard Brownian motion. Then $R$ admits the following Lamperti representation:

$$R_t = x \exp B^{(\nu)}_{\tau^{(\nu)}(tx^{-1/2})}, \qquad t \geq 0,$$

where $B^{(\nu)}$ is a Brownian motion with drift $\nu = \mu/2 - 1$, that is, $B^{(\nu)} = (B_t + \nu t, t \geq 0)$, and $\tau_t^{(\nu)} = \inf\{s : \int_0^s \exp B_s^{(\nu)} > t\}$. Since $\limsup_t B_t^{(\nu)} = +\infty$ if and only if $\nu \geq 0$, the cases which are treated in this paper concern the Bessel process with dimension $\mu \geq 2$, raised to any positive power. (Note that the only case where the process $R$ does not drift to $+\infty$ is when $\mu = 2$.) Then it is well known that conditions of Theorems 1 and 2 are satisfied and that the limiting process which is defined in Theorem 1 is the unique strong solution of (4.1), for $x = 0$.

The class of stable Lévy processes conditioned to stay positive provides another interesting example. Let $Y$ be any such process with index $\alpha \in (0, 2]$, and let $P_x$ be the law of $Y + x$. The real function $x \mapsto x^{\alpha\rho}$ defined on $[0, \infty)$, where $\rho = P_0(Y_1 \leq 0)$, is harmonic for the semigroup of the process $Y$ killed at its first entrance time into the negative half-line. Let $\tau = \inf\{t : Y_t \leq 0\}$ be this time, and call $(\mathcal{F}_t)$ the usual filtration on $\mathcal{D}$, then the $h$-process whose law is defined as follows:

$$\mathbb{P}_x(\Lambda) = x^{-\alpha\rho} E_x(Y_t^{\alpha\rho} \mathbb{1}_\Lambda \mathbb{1}_{\{t < \tau\}}), \qquad \Lambda \in \mathcal{F}_t, t \geq 0, x > 0,$$

is a strong Markov process which is called the Lévy process $Y$ conditioned to stay positive. A more intuitive way to define this process is to first condition $Y$ to stay positive on the time interval $[0, s]$ and then to let $s$ go to $+\infty$:

$$\mathbb{P}_x(\Lambda) = \lim_{s \uparrow \infty} P_x(\Lambda | \tau > s), \qquad \Lambda \in \mathcal{F}_t, t \geq 0, x > 0.$$

It is clear that $(X, \mathbb{P}_x)$ is a pssMp. We emphasize that when $Y$ is the standard Brownian motion, the process $(X, \mathbb{P}_x)$ corresponds to the three-dimensional Bessel process. We refer to [6] for a more complete study of Lévy processes conditioned to stay positive. In particular, in [6], it is shown that $(X, \mathbb{P}_x)$ satisfies $\lim_{t \uparrow \infty} X_t = +\infty$, $\mathbb{P}_x$-a.s. Moreover, $(X, \mathbb{P}_x)$ reaches its overall minimum only once and $(X, \mathbb{P}_x)$ converges weakly on $\mathcal{D}$, as $x$ tends to 0, toward the law of the post-minimum process $(X_{m+t}, t \geq 0)$, $m = \inf\{t : X_t = \inf_{s \geq 0} X_s\}$, see also [7]. As a pssMp, $(X, \mathbb{P}_x)$ admits a Lamperti representation. The law of the underlying Lévy process $\xi$ in this representation has been computed in [5] in terms of the law of $Y$. Suppose that $\alpha \in [1, 2)$. Let $\Phi$ be the characteristic exponent of $Y$, that is, $\Phi(u) = -\log E(\exp iuY_1)$, then it is well known that $\Phi$ admits the Lévy Khintchnie decomposition

$$\Phi(u) = imu + \int_{(-\infty,\infty)} (1 - e^{iux} + iux\mathbb{1}_{\{|x|\leq 1\}})\mu(dx),$$

where $m \in \mathbb{R}$ and

$$\mu(dx) = (c_+ x^{-(\alpha+1)}\mathbb{1}_{\{x>0\}} + c_- |x|^{-(\alpha+1)}\mathbb{1}_{\{x<0\}})\, dx,$$

$c_+$ and $c_-$ being positive constants. According to [5], the characteristic exponent of $\xi$, introduced in Section 2.1, admits the Lévy Khintchine decomposition

$$\psi(u) = iau + \int_{(-\infty,\infty)} (1 - e^{iux} + iux\mathbb{1}_{\{|x|\leq 1\}})\pi(dx),$$

with $a \in \mathbb{R}$ and with Lévy measure

$$\pi(dx) = \left(\frac{c_+ e^{(\alpha\rho+1)x}\mathbb{1}_{\{x>0\}}}{(e^x-1)^{\alpha+1}} + \frac{c_- e^{(\alpha\rho+1)x}\mathbb{1}_{\{x<0\}}}{|e^x-1|^{\alpha+1}}\right) dx.$$

Here again, it can be checked that the conditions of Theorems 1 and 2 are satisfied.

**Acknowledgments.** We are indebted to Ron Doney for the private communication [9] and many helpful discussions on overshoots of Lévy processes. The present work has been done during the visit of L. Chaumont to the University of Mexico UNAM and the one of M. E. Caballero to the University of Paris X. We thank both of these universities for their support.

Instituto de Matemáticas
Universidad Nacional Autónoma de México
Mexico 04510 DF
E-mail: emilia@servidor.unam.mx

Laboratoire de Probabilités
  et Modèles Aléatoires
Université Pierre et Marie Curie
4 Place Jussieu 75252
Paris Cedex 05
France
E-mail: chaumont@ccr.jussieu.fr
URL: www.proba.jussieu.fr/pageperso/chaumont